\documentclass[12pt,a4paper]{amsart}

\usepackage[all]{xy}

\newtheorem{theorem}{Theorem}

\newtheorem{lemma}[theorem]{Lemma}

\theoremstyle{remark}

\newtheorem*{problem*}{Problem}
\newtheorem*{remark*}{Remark}
\newtheorem*{convention*}{Convention}
\newtheorem*{notation*}{Notation}
\newtheorem*{examples*}{Examples}
\newtheorem*{example*}{Example}
\newtheorem*{warning*}{Warning}

\def\N{{\mathbb N}}

\def\R{{\mathbb R}}
\def\T{{\mathbb T}}
\def\C{{\mathbb C}}
\def\Z{{\mathbb Z}}

\newcommand{\range}{\operatorname{range}}

\begin{document}

\title[From filters to wavelets via direct limits]{From filters to wavelets via direct limits}

\author[Nadia S. Larsen]{Nadia S. Larsen}

\address{Mathematics Institute, University of Oslo, Blindern, NO-0316 Oslo, Norway}
\email{nadiasl@math.uio.no}

\author[Iain Raeburn]{Iain Raeburn}

\address{School  of Mathematical and Physical Sciences, University of
Newcastle, NSW 2308, Australia}

\email{iain.raeburn@newcastle.edu.au}

\begin{abstract}
We present a new proof of a theorem of Mallat which describes a construction of wavelets starting from a quadrature mirror filter. Our main innovation is to show how the scaling function associated to the filter can be used to identify a certain direct limit of Hilbert spaces with $L^2(\R)$ in such a way that one can immediately identify the wavelet basis. Our arguments also use a pair of isometries introduced by Bratteli and Jorgensen, and exploit the geometry inherent in the Cuntz relations satisfied by these isometries.
\end{abstract}

\date{18 August 2005}

\thanks{The authors thank Ola Bratteli and Palle Jorgensen for interesting conversations in Oslo in June 2004 (Palle's slant on the ideas we discussed is incorporated in~\cite{dutjor}). Iain Raeburn also wishes to thank Larry Baggett, Astrid an Huef, Kathy Merrill, Judy Packer and Arlan Ramsay for helpful conversations about recent extensions of Mallat's theorem (as in, for example, \cite{bjmp}).} 

\thanks{This research was supported by the Australian Research Council, through the ARC Centre for Complex Dynamic Systems and Control, and the Norwegian Research Council.}

\maketitle

A \emph{wavelet} is a function $\psi\in L^2(\R)$ such that
\[
\{\psi_{j,k}:x\mapsto 2^{-j/2}\psi(2^{-j}x-k):j,k\in Z\}
\]
is an orthonormal basis for $L^2(\R)$. There are, remarkably, many different wavelets, and they have proved to be enormously useful in both theory and applications. So there has been a great deal of interest in methods of constructing wavelets. One famous construction of Mallat~\cite{mal} starts from a \emph{quadrature mirror filter}: a function $m_0:\T\to \C$ such that 
\[
|m_0(z)|^2+|m_0(-z)|^2=1 \ \mbox{for every $z\in \T$.}
\]
Our goal here is to present a new proof of Mallat's theorem based on the concept of a direct limit.

Mallat proved his theorem in two stages. From the filter $m_0$ he built a \emph{multiresolution analysis}, in which a central role is played by a \emph{scaling function} $\phi\in L^2(\R)$ satisfying $\phi(2x)=m_0(e^{2\pi ix})\phi(x)$ \cite[Theorem~2]{mal}. He then used what he described as a ``by now classic'' algorithm to generate the wavelet \cite[\S4]{mal}. Mallat's construction has since been refined and discussed in several books. For example, \cite[\S5.3--4]{fraz} contains a relatively elementary proof of his theorem, in which some of the analysis has been simplified but the overall strategy is that of Mallat. In our proof, the scaling function still plays a central role: we use it to identify a certain direct limit with $L^2(\R)$, and the existence of the wavelet then follows almost immediately from the geometry implicit in some operator-theoretic equations called the \emph{Cuntz relations}. The analytic content of our proof is much the same as that in the standard sources, and we refer to them for details, but our organisation seems to be quite different.  

Our arguments seem to be more natural in the Fourier or frequency domain, so we work there throughout, and our construction yields the Fourier transform of the wavelet. One effect of working in the Fourier domain is that the scaling equation (Equation~\eqref{scaling1} below) involves multiplication rather than convolution. In the interests of clarity, we shall consider only the classical (dyadic) wavelets, but we are optimistic that our approach will also shed light in other situations where wavelet bases are used.

\medskip

Throughout, $m_0$ will be a quadrature mirror filter such that $m_0$ is smooth at $1$, $m_0(1)=1$ and $m_0(z)\not=0$ for $z$ in the right half-circle. We define $m_1:\T\to \C$ by \[m_1(z):=z\overline{m_0(-z)}.\]
We now define two operators $S_0$ and $S_1$ on $L^2(\T)$ by
\begin{equation}\label{defSi}
(S_if)(z):=2^{1/2}m_i(z)f(z^2).
\end{equation}
Our starting point is the following observation of Bratteli and Jorgensen \cite{bj}:

\begin{lemma}\label{bjlem}
The operators $S_i$ satisfy $S_0^*S_0=1=S_1^*S_1$ and $S_0S_0^*+S_1S_1^*=1$.
\end{lemma} 

To prove this, first verify that the adjoints $S_i^*$ are given by
\[
(S_i^*f)(e^{2\pi ix})=2^{-1/2}\big(\,\overline{m_i(e^{\pi ix})}f(e^{\pi ix})+\overline{m_i(-e^{\pi ix})}f(-e^{\pi ix})\big),
\]
and then compute $S_i^*S_i$ and $S_0S_0^*+S_1S_1^*$. 

The formal computation in Lemma~\ref{bjlem} has some very interesting geometric consequences. The relations $S_i^*S_i=1$ say that the operators $S_i$ are isometries of $L^2(\T)$ into itself, and imply that the operators $S_iS_i^*$ are the orthogonal projections onto the ranges of the $S_i$ (which are automatically closed because the $S_i$ are isometries). Since the sum of two projections is a projection only when their ranges are orthogonal, the \emph{Cuntz relation} $S_0S_0^*+S_1S_1^*=1$ implies that the ranges $S_iL^2(\T)$ of the $S_i$ are orthogonal complements of each other, so that 
\begin{equation}\label{decompL2}
L^2(\T)=S_0L^2(\T)\oplus S_1L^2(\T).
\end{equation}

We now recall a construction from algebra. Suppose that we have Hilbert spaces $H_n$ and isometries $T_n:H_n\to H_{n+1}$ for all $n\in \N$. The \emph{direct limit} $(H_\infty,U_n)$ consists of a Hilbert space $H_\infty$ and isometries $U_n:H_n\to H_\infty$ which satisfy $U_{n+1}\circ T_n=U_n$ and which have the following universal property: for every family of isometries $\{R_n\}$ of $H_n$ into a Hilbert space $K$ such that $R_{n+1}\circ T_n=R_n$, there is a unique isometry $R_\infty:H_\infty\to K$ such that $R_\infty\circ T_n=R_n$ for every $n$. We illustrate this universal property with the diagram:
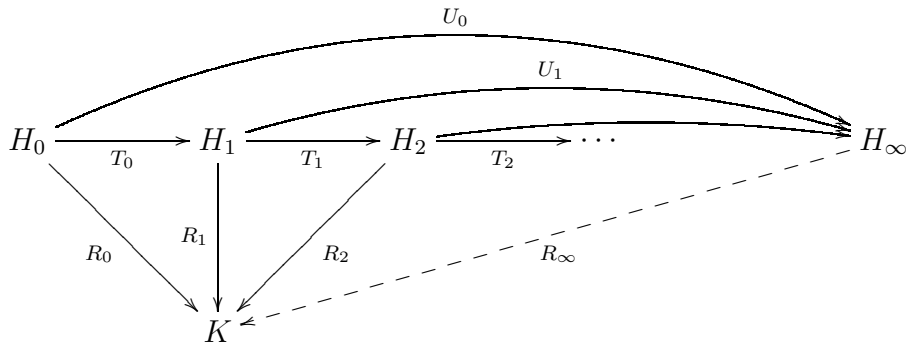
\[\xygraph{
{H_0}="v0":[rr]{H_1}="v1"_{T_0}:[rr]{H_2}="v2"_{T_1}:[rr]{\cdots}="v3"_{T_2}:@{}[rrr]
{H_\infty}="v4"
"v0":@/^{40pt}/^{U_0}"v4""v1":@/^{20pt}/"v4"^{U_1}"v2":@/^{7pt}/"v4"^{}
"v1":[dd]{K}="w1"_{R_1}"v0":"w1"_{R_0}"v2":"w1"^{R_2}"v4":@{-->}"w1"^{R_\infty}
}\]
The uniqueness implies that $H_\infty=\overline{\bigcup_{n=1}^\infty U_nH_n}$, since otherwise we could define $R_\infty$ arbitrarily on $(\bigcup U_nH_n)^\perp$. 

It is not hard to see that the direct limit exists. Indeed, identifying each $h\in H_n$ with all its images $T_{m-1}T_{m-2}\cdots T_nh\in H_m$ defines an equivalence relation $\sim$ on the disjoint union $\bigsqcup H_n$, and because the $T_n$ are isometries and hence inner-product preserving, the set $H':=(\bigsqcup H_n)/\!\!\sim$ of equivalence classes is naturally an inner-product space; completing $H'$ gives a Hilbert space $H_\infty$, and the maps $U_n$ which send elements of $H_n$ to their class in $H'\subset H_\infty$ have the required properties. However, the important point is that the construction does not matter, since the universal property identifies the direct limit up to isomorphism: to identify $H_\infty$ with a Hilbert space $K$, for example, we just need to find isometries $R_n:H_n\to K$ as above such that $K=\overline{\bigcup_{n=1}^\infty R_nH_n}$, and then $R_\infty$ is an isomorphism of $H_\infty$ onto $K$.

If we start with a single isometry $S$ on a Hilbert space $H$, we can take the direct limit $(H_\infty, U_n)$ of the system in which every $H_n$ is $H$ and every $T_n$ is $S$. Now consider the diagram 
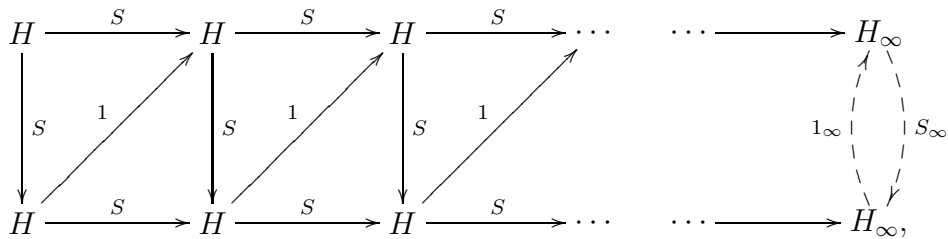
\[\xygraph{
{H}="v0":[rr]{H}="v1"^{S}:[rr]{H}="v2"^{S}:[rr]{\cdots}="v3"^{S}:@{}[r]{\cdots}="v_3":[rr]
{H_\infty}="v4"
"v0":[dd]{H}="w0"^S
"w0":[rr]{H}="w1"^{S}:[rr]{H}="w2"^{S}:[rr]{\cdots}="w3"^{S}:@{}[r]{\cdots}="w_3":[rr]
{H_\infty,}="w4"
"v1":"w1"^S"v2":"w2"^S"v4":@{-->}@/^{10pt}/"w4"^{S_\infty}
"w0":"v1"^1"w1":"v2"^1"w2":"v3"^1"w4":@{-->}@/^{10pt}/"v4"^{1_\infty}
}\]
where the horizontal arrows going into $H_\infty$ are there to remind us that we have isometries $U_n$ of every copy $H_n$ of $H$ into $H_\infty$. Applying the universal property of the top row to the downward arrows (or more properly, to $R_n:=U_n\circ S$) gives an isometry $S_\infty$ on $H_\infty$ which is characterised by $S_\infty(U_nh)=U_n(Sh)$. Similarly, applying the universal property of the bottom row to the NE arrows gives an isometry $1_\infty$ which is characterised by $1_\infty(U_{n}h)=U_{n+1}h$. Then for every $n$ and every $h\in H$ we have
\begin{align*}
(1_\infty S_\infty)(U_nh)&=1_\infty(S_\infty U_nh)=1_\infty(U_n(Sh))=U_{n+1}(Sh)=U_nh\ \mbox{ and}\\
(S_\infty 1_\infty)(U_nh)&=S_\infty(U_{n+1}h)=U_{n+1}(Sh)=U_nh,
\end{align*}
which implies that $1_\infty$ is an inverse for $S_\infty$. This process of passing to the direct limit, therefore, turns the isometry $S$ into a unitary $S_\infty$. Since 
\begin{equation}\label{Sinftyscales}
S_\infty(U_{n+1}h)=U_{n+1}(Sh)=U_nh,
\end{equation} 
this unitary is an isomorphism of the copy $U_{n+1}H$ of $H_{n+1}=H$ onto the copy $U_nH$ of $H_n=H$.

Applying the process described in the previous paragraph to the isometry $S_0$ on $L^2(\T)$ defined in \eqref{defSi} gives a Hilbert space $H_\infty$ and a unitary operator $S_\infty$ on $H_\infty$. Our next task is to identify the direct limit $H_\infty$ with $L^2(\R)$. Mallat proved that, under our hypotheses on $m_0$, there is a \emph{scaling function}\footnote{In the literature, it is the inverse Fourier transform $\check \phi$ of this function $\phi$ which is usually called a scaling function.} $\phi\in L^2(\R)$ of norm $1$ such that
\begin{gather}
\label{scaling1}\phi(2x)=m_0(e^{2\pi ix})\phi(x)\ \mbox{ and}\\
\label{scaling2}\sum_{k\in \Z}|\phi(x+k)|^2=1
\end{gather}
for every $x\in \R$. Indeed, he proved that the infinite product 
\[
\phi(x)=\prod_{n=1}^\infty m_0(\exp(2\pi i2^{-n}x))
\]
has the required properties (see \cite[pages 76--77]{mal} or \cite[\S5.4]{fraz}). For example, if $m_0$ is the characteristic function of the right half-circle, then $\phi=\chi_{[-1/2,1/2]}$ is a scaling function. (In \cite[page~225]{hol} this is proved under the milder hypothesis that $m_0$ satisfies \emph{Cohen's condition}, which is a necessary and sufficient condition for the existence of such a function $\phi$.) 

We now fix $n\in \Z$, and define $R_n:H_n=L^2(\T)\to L^2(\R)$ by
\[
(R_nf)(x)=2^{-n/2}f(\exp(2\pi i2^{-n}x))\phi(2^{-n}x).
\]
Each $R_n$ is an isometry: indeed, a change of variables and an application of the monotone convergence theorem shows that
\begin{align*}
\|R_nf\|_{L^2(\R)}^2&=\int_{\R} 2^{-n}|f(\exp(2\pi i2^{-n}x))\phi(2^{-n}x)|^2\,dx\\
&=\int_{\R} |f(e^{2\pi ix})\phi(x)|^2\,dx\\
&=\sum_{n\in \Z}\int_0^1|f(e^{2\pi it})\phi(t+n)|^2\,dt\\
&=\int_0^1|f(e^{2\pi it})|^2\Big(\sum_{n\in\Z}|\phi(t+n)|^2\Big)\,dt,
\end{align*}
which by \eqref{scaling2} is precisely the norm of $f$ in $L^2(\T)$. The identity \eqref{scaling1} implies that the isometries $R_n$ are compatible with the isometries $S_0:H_n=L^2(\T)\to H_{n+1}=L^2(\T)$ in the direct system:
\begin{align*}
R_{n+1}(S_0f)(x)&=2^{-(n+1)/2}2^{1/2}m_0(\exp(2\pi i2^{-(n+1)}x))f(\exp(2\pi i2^{-(n+1)}x)^2)\phi(2^{-(n+1)}x)\\
&=2^{-n/2}f(\exp(4\pi i2^{-(n+1)}x))\phi(2(2^{-(n+1)}x))\\
&=(R_nf)(x).
\end{align*} 
Thus the universal property of the direct limit gives an isometry $R_\infty:H_\infty\to L^2(\R)$. Notice that if we set $V_{n}:=\range R_n=R_\infty(U_nL^2(\T))$, then the last calculation shows that $V_{n}\subset V_{n+1}$.

\begin{lemma}\label{idSinfty}
The isometry $R_\infty$ intertwines the unitary $S_\infty$ on $H_\infty$ and the dilation operator $D$ on $L^2(\R)$ defined by $(D\xi)(x)=2^{1/2}\xi(2x)$.
\end{lemma}

\begin{proof}
We let $f\in H_{n+1}=L^2(\T)$, and compute:
\begin{align*}
D(R_\infty(U_{n+1}f))(x)&=D(R_{n+1}f)(x)\\
&=2^{1/2}(R_{n+1}f)(2x)\\
&=2^{1/2}2^{-(n+1)/2}f(\exp(2\pi i2^{-(n+1)}2x))\phi(2^{-(n+1)}2x)\\
&=(R_nf)(x)\\
&=R_\infty(U_nf)(x).
\end{align*}
Since we know from \eqref{Sinftyscales} that $S_\infty\circ U_{n+1}=U_n$, this implies that $D\circ R_\infty$ and $R_\infty\circ S_\infty$ agree on the range of every $U_{n+1}$, and hence on $H_\infty$.
\end{proof}

Since $S_\infty$ is an isomorphism of $U_{n+1}H$ onto $U_nH$, it follows from Lemma~\ref{idSinfty} that\footnote{This formula looks slightly different from the usual dilation property of a multiresolution analysis because we are working in the frequency domain; in \cite{mal} and \cite{fraz}, for example, the space denoted by $V_n$ is the inverse Fourier transform of our $V_n$.} $DV_{n+1}=V_{n}$ for every $n\in \Z$. The following lemma is proved in \cite[pages~78--79]{mal} and \cite[Lemmas~5.47 and 5.48]{fraz}. (A different argument which proves the analogous property of the inverse Fourier transforms is given in Propositions~5.3.1 and 5.3.2 of~\cite{dau}.)

\begin{lemma}\label{Visbig}
We have $\bigcap_{n\in \Z}V_n=\{0\}$ and $\overline{\bigcup_{n\in \Z}V_n}=L^2(\R)$.
\end{lemma}

We write $W_n$ for the complement $V_{n+1}\ominus V_n$ of $V_n$ in $V_{n+1}$.  Then the subspaces $W_n$ are mutually orthogonal, and it follows from Lemma~\ref{Visbig} that $L^2(\R)$ decomposes as the direct sum $\bigoplus_{n\in \Z}W_n$. Since $D^{-1}V_{n-1}=V_n$ for every $n$, we have $DW_{n+1}=W_n$ for every $n$. Thus to find an orthonormal basis for $L^2(\R)$, it suffices to find an orthonormal basis for one $W_n$, and then this together with all its dilates will be an orthonormal basis for $L^2(\R)$.

So we seek an orthonormal basis for $W_0$. With $H=L^2(\T)$, we have
\[
W_0=V_{1}\ominus V_0=R_{1}H\ominus R_0H=R_{1}H\ominus R_{1}(S_0H)=R_{1}(H\ominus S_0H).
\]
At this point we recall from \eqref{decompL2} that the complement of $S_0H$ is the range $S_1H$ of the other isometry $S_1$. Since $R_{1}S_1$ is an isometry, it maps the usual orthonormal basis $\{e_k:z\mapsto z^k:k\in \Z\}$ for $L^2(\T)$ into an orthonormal basis for $W_0$. We deduce that the functions 
\[
\psi_k(x)=R_{1}S_1(e_{-k})(x)=2^{-1/2}\big(2^{1/2}m_1(\exp(2\pi i2^{-1}x))e^{-2\pi ikx}\big)\phi(2^{-1}x)
\]
form an orthonormal basis for $W_0$. We set 
\[
\psi(x):=m_1(e^{\pi ix})\phi(2^{-1}x),
\]
so that the basis elements take the form
\[
\psi_k(x)=e^{-2\pi ikx}\psi(x).
\]
If we now define
\[
\psi_{j,k}(x)=(D^{j}\psi_k)(x)=2^{j/2}\exp({-2\pi ik2^{j}x})\psi(2^{j}x),
\]
then $\{\psi_{j,k}:j,k\in \Z\}$ is an orthonormal basis for $L^2(\R)$. 

Since the inverse Fourier transform intertwines $D$ and $D^{-1}$, and intertwines multiplication by $e^{2\pi ikx}$ and the translation operator which takes $\xi$ to $\xi(\cdot+k)$, the functions
\[
\check\psi_{j,k}(x)=2^{-j/2}\check\psi(2^{-j}x-k)
\]
also form an orthonormal basis for $L^2(\R)$. In other words, $\check\psi$ is a wavelet, and we have proved Mallat's theorem:

\begin{theorem}
Suppose that $m_0$ is a quadrature mirror filter such that $m_0$ is smooth at $1$, $m_0(1)=1$ and $m_0(z)\not=0$ for $z$ in the right half-circle, and let $\phi$ be a function satisfying the scaling conditions \eqref{scaling1} and \eqref{scaling2}. Define $\psi:\R\to \C$ by
\[
\psi(x)=e^{\pi ix}\,\overline{m_0(-e^{\pi ix})}\,\phi(x/2).
\]
Then the inverse Fourier transform $\check \psi$ is a wavelet.
\end{theorem}

\end{document}